\def\ds{\displaystyle}
\def\RR{{\mathbb R}}

\def\epsilon{\varepsilon}
\documentstyle[12pt]{amsart}
\title{Variations on an inequality from IMO'2001}
\author{Oleg Mushkarov \& Nikolai Nikolov}
\address
{Institute of Mathematics and Informatics\\
Bulgarian Academy of Sciences\\1113 Sofia, Bulgaria}
\email{muskarov@@math.bas.bg, nik@@math.bas.bg}
\begin{document}
\begin{abstract}
Some extensions of an inequality from IMO'2001 are proven by means of the
Lagrange multiplier criterion.
\end{abstract}
\maketitle

\section{Introduction}

This paper is a continuation of \cite{1,2}, where some natural generalizations
of Problem 2 from IMO'2001 have been proved. Oir aim here is to consider
some other extensions of the same problem which states:

{\it Prove that 
$$\frac{a}{\sqrt{a^2+8bc}}+\frac{b}{\sqrt{b^2+8ac}}+\frac{c}{\sqrt{c^2+8ab}}
\geq 1,$$
where $a,b$ and $c$ are arbitrary positive numbers.} 

Many different proofs of this inequality were given during the Olym-\\piad
and it was  also shown by the first author that the following more general 
inequality holds:
\begin{equation}
\frac{a}{\sqrt{a^2+\lambda bc}}+\frac{b}{\sqrt{b^2+\lambda ac}}+
\frac{c}{\sqrt{c^2+\lambda ab}}\geq \frac{3}{\sqrt{1+\lambda}}
\end{equation}
for arbitrary $a,b,c>0$ and $\lambda \geq 8.$ It is easy to see that the 
latter inequality is not true for $0<\lambda <8$. Moreover, it can be shown that
in this case the infimum of the function in the left-hand side of (1) (when 
$a,b$ and $c$ run over all positive numbers) is equal to 1. This fenomenon 
led us to consider the following general problem:

{\it Find the infimum and the supremum of the function
$$F_{\alpha}(x_1,x_2,\dots , x_n)=\sum_{i=1}^{n} \frac{1}{(1+x_i)^\alpha}$$
on the set 
$$H_\lambda =\{(x_1,x_2,\dots, x_n)\in\RR^n| x_1x_2\dots x_n=\lambda^n,
x_1,x_2, \dots , x_n>0\},$$
where $\lambda >0$ and $\alpha$ are given real constants.} 

\section{The infimum of $F_\alpha$}

We shall find the infimum of the function $F_\alpha$ on the set $H_\alpha$ by
means of the well-known Lagrange multiplier criterion. The next proposition has
been proved in \cite{2}, but we include it here to make the paper self-contained.

\

{\bf Proposition 1.} {\it For any $\alpha\in (0,1]$ we have 
$$\inf_{H_\lambda}F_\alpha=\min(1,\frac{n}{(1+\lambda)^\alpha}).$$}

\

{\bf Proof.} Suppose first that $d:=\ds\inf_{H_\lambda}F_\alpha$ is not attained
at a point of $H_\lambda$. Then, 
$d=F_\alpha=\lim_{k\to \infty}F_\alpha(x_1^{(k)},\dots,x_n^{(k)}),$ where,
for example, $\ds\lim_{k\to\infty}x_n^{(k)}=0$ or $+\infty.$ Hence, for example,
$\ds\lim_{k\to\infty}x_1^{(k)}=+\infty$ or $0$ and in both cases we see that 
$d\geq 1.$ Note that if $\ds\lim_{k\to\infty}x_s^{(k)}=+\infty$ 
for $s=1,2\dots,n-1$ and $\ds\lim_{k\to\infty}x_n^{(k)}=0,$ then 
$\ds\lim_{k\to\infty}F_\alpha(x_1^{(k)},\dots,x_n^{(k)})=1.$ Now, let 
$d$ is attained at a point of $H_\lambda$. Consider the function 
$$F(x_1,x_2,\dots ,x_n)=F_\alpha(x_1,x_2,\dots ,x_n)+
\mu(x_1x_2\dots x_n-\lambda^n).$$
Then the Lagrange multiplier criterion says that $d$ is attained at a point
$(x_1,x_2,\dots,x_n)\in H_\lambda$ such that 
$$\frac{\partial F}{\partial x_i}=-\frac{\alpha}{(1+x_i)^{\alpha+1}}+
\frac{\mu x_1\dots x_n}{x_i}=0,$$
i.e., when 
\begin{equation}
\frac{x_i}{(1+x_i)^{\alpha+1}}=\frac{x_j}{(1+x_j)^{\alpha+1}}, 1\leq i,j\leq n.
\end{equation}
Consider the function $\ds g(x)=\frac{x}{(1+x)^{\alpha+1}}.$ Then, 
$\ds g'(x)=\frac{1-\alpha x}{(1+x)^{\alpha+2}},$ and, therefore, $g(x)$ takes
each its value at most twice. Hence (2) shows that $x_1=\cdots =x_k=x$ and
$x_{k+1}=\cdots=x_n=y$ for some $1\leq k\leq n.$ If $k=n,$ then
$x_1=x_2=\cdots=x_n=\lambda$ and
$\ds F_\alpha(x_1,x_2,\dots,x_n)=\frac{n}{(1+\alpha)^\lambda}.$
If $k<n,$ then
$$F_\alpha(x_1,x_2,\dots,x_n)=\frac{k}{(1+x)^\alpha}+\frac{n-k}{(1+y)^\alpha}
\geq\frac{1}{(1+x)^\alpha}+\frac{1}{(1+y)^\alpha}.$$
To prove Proposition 1 it is sufficient to show that 
\begin{equation}
\frac{1}{(1+x)^\alpha}+\frac{1}{(1+y)^\alpha}>1
\end{equation}
provided
\begin{equation}
\frac{x}{(1+x)^{\alpha+1}}=\frac{y}{(1+y)^{\alpha+1}},x\neq y.
\end{equation}
Set $\ds\beta=\frac{1}{\alpha}\geq 1,z=(1+x)^\alpha$ and $t=(1+y)^\alpha.$ 
Then (3) and (4) can be written respectively as $z+t>zt$ and
$\ds(zt)^\beta=\frac{z^{\beta+1}-t^{\beta+1}}{z-t}.$ So, we have to prove that 
\begin{equation}
(z+t)^\beta\geq\frac{z^{\beta+1}-t^{\beta+1}}{z-t}.
\end{equation}
Assume that $z<t$ and set $\ds u=\frac{z}{t}<1.$ Applying Bernoulli's inequality
twice we obtain $\ds(1+u)^\beta\ge 1+\beta u>\frac{1-u^{\beta+1}}{1-u}$ which
is just the inequality (5).\qed

\

The next example shows that for a given $\alpha>1$ a result similar to 
Proposition 1 could be expected only for sufficiently large $n.$ 

\

{\bf Example 1.} {\it Let $\alpha=2$ and $n=2.$ Then the function
$F_2(x_1,x_2)$ attains minimum on $H_\lambda$ given by
\begin{equation}
\min_{H_\lambda}F_2=\left\{\begin{array}{cc}
\ds\frac{2}{(1+\lambda)^2}&\mbox{if }\lambda\geq\frac{1}{2}\\
\ds\frac{1-2\lambda^2}{(1-\lambda^2)^2}&\mbox{if }0<\lambda\leq\frac{1}{2}.
\end{array}\right.
\end{equation}}

\

{\bf Proof.} To prove (6) we proceed as in the proof of Proposition 1. First
note that if $x_1\to 0$ or $+\infty,$ then  $x_2\to+\infty$ or $0$ and, in both
cases, $F_2(x_1,x_2)\to 1.$ Consider the points $(x_1,x_2)\in H_\lambda$ such that
\begin{equation}
\frac{x_1}{(1+x_1)^3}=\frac{x_2}{(1+x_2)^3}.
\end{equation}
If $x_1=x_2=\lambda,$ then $\ds F_2(x_1,x_2)=\frac{2}{(1+\lambda)^2}.$ If 
$x_1\neq x_2,$ then (7) is equivalent to $\ds x_1+x_2=\frac{1}{\lambda^2}-3.$
This together with $x_1x_2=\lambda^2$ implies that 
$\ds\frac{1}{\lambda^2}-3\geq 2\lambda,$ i.e., $\ds\lambda<\frac{1}{2}$ and 
$\ds F_2(x_1,x_2)=\frac{1-2\lambda^2}{(1-\lambda^2)^2}.$ Hence (6) follows 
from the inequalities $\ds\frac{1-2\lambda^2}{(1-\lambda^2)^2}<1$ and 
$\ds\frac{1-2\lambda^2}{(1-\lambda^2)^2}<\frac{2}{(1+\lambda)^2}$ for any 
$\lambda >0$, and $\ds\frac{2}{(1+\lambda)^2}<1$ for $\lambda\geq\frac{1}{2}.$
\qed

\

The next proposition gives a partial result in the case $\alpha>1.$

\

{\bf Proposition 2.} {\it For any $\alpha >1$ and any integer $n\ge\alpha+1$ 
we have 
$$\inf_{H_\lambda}F_\alpha=\min(1,\frac{n}{(1+\lambda)^\alpha}).$$}

\

{\bf Proof.} Proceeding as in the proof of Proposition 1 it is sufficient to
prove that 
$$(1+(n-1)u)^\beta>\frac{1-u^{\beta+1}}{1-u}$$
for $\ds\beta=\frac{1}{\alpha}<1$ and $0<u<1.$ Since $n-1\geq\alpha$ we have
$\ds 1+(n-1)u\geq 1+\frac{n}{\beta}$ and it is enough to show that 
\begin{equation}
(1+\frac{u}{\beta})^\beta>\frac{1-u^{\beta+1}}{1-u}
\end{equation}
for $\beta,u\in(0,1)$. Consider the function
$$f(x)=(1-x)(1+\frac{x}{\beta})^\beta+x^{\beta+1}-1\hbox{ for }x\in[0,1].$$
Since  
$$f'(x)=\frac{(1+\beta)x}{\beta}(\beta x^{\beta-1}-(1+\frac{x}{\beta})^{\beta-1})$$
the equation $f'(x)=0$ has a unique real root   
$\ds x_0=(\beta^{\frac{1}{\beta-1}}-\beta^{-1})^{-1}.$
On the other hand, since $f(0)=f(1)=0$ and $\beta-1<0,$ it follows that
$x_0\in(0,1),\ f'(x)>0$ for $x\in(0,x_0)$ and $f'(x)<0$ for $x\in(x_0,1).$
Hence $f(x)>0$ for $x\in (0,1)$ and the inequality (8) is proved.\qed 

\

{\bf Remark 1.} As Example 1 suggests, if $\alpha>1$ and $n<\alpha+1,$ then a
result similar to Proposition 2 is not true. The authors do not know the value
of $\ds\inf_{H_\lambda}F_\alpha$ for such $\alpha$ and $n.$

\

To complete this section it remains to consider the case $\alpha <0$.

\

{\bf Proposition 3.} {\it For any $\alpha<0$ the function
$F_\alpha(x_1,\dots, x_n)$ attains minimum on $H_\lambda$ given by 
$$\min_{H_\lambda}F_\alpha =\frac{n}{(1+\lambda)^\alpha}.$$}

\

{\bf Proof.} We may proceed as in the proof of Proposition 1 but in this case 
the statement follows directly from the fact that the function 
$\ds f(x)=\frac{1}{(1+e^x)^\alpha}$ is convex for $\alpha<0$ since
$f''(x)>0.$\qed

\section{The supremum of $F_\alpha$}

The results obtained in this section are dual analogs of that in Section 2.

\

{\bf Proposition 4.} {\it For any $\alpha\geq 1$ we have 
$$\sup_{H_\lambda}F_\alpha=\max(n-1,\frac{n}{(1+\lambda)^\alpha}).$$}

\

{\bf Proof.} We proceed as in the proof of Proposition 1. If    
$\ds\sup_{H_\lambda}F_\alpha$ is not attained at a point of $H_\lambda$ then
we may assume that $x_n\to+\infty$ and obviously we have
$\ds\sup_{H_\lambda}F_\alpha\leq n-1$. Note also that if $x_1\to0,\dots, 
x_{n-1}\to 0$ and $x_n\to+\infty,$ then $F_\alpha(x_1,\dots ,x_n)\to n-1.$

Next consider the case when $\ds\sup_{H_\lambda}F_\alpha$ is attained at
a point of $H_\lambda$ such that $x_1=\dots=x_k=x$ and
$x_{k+1}=\dots=x_n=y$. If $x=y,$ then $x_1=\dots =x_n=\lambda$ and
$\ds F_\alpha(x_1,\dots,x_n)=\frac{n}{(1+\lambda)^\alpha}.$ If
$x\neq y,$ then $k<n$ and $\ds F_\alpha(x_1,\dots,x_n)=\frac{k}{(1+x)^\alpha}
+\frac{n-k}{(1+y)^\alpha}$. So, it is enough to prove that if
$\ds\frac{x}{(1+x)^{\alpha+1}}=\frac{y}{(1+y)^{\alpha+1}}$ and $x<y,$ then
$\ds\frac{n-1}{(1+x)^{\alpha}}+\frac{1}{(1+y)^{\alpha}}<n-1$. But this follows
from the inequality $\ds\frac{1}{(1+x)^{\alpha}}+\frac{1}{(1+y)^{\alpha}}<1,$
which can be proved by using Bernoulli's inequality for $\ds\beta=\frac{1}{\alpha}<1$
as in the proof of Proposition 1.\hfill$\Box$

\

The next example is dual to Example 1.

\

{\bf Example 2.} {\it Let $\ds\alpha=\frac{1}{2}$ and $n=2$. Then the function
$F_{\frac{1}{2}}(x_1,x_2)$ attains maximum on $H_\lambda$ given by
\begin{equation}
\max_{H_\lambda}F_{\frac{1}{2}}=\left\{\begin{array}{cc}
\ds\frac{\lambda}{\sqrt{\lambda^2-1}}&\mbox{if }\lambda>2\\
\ds\frac{2}{\sqrt{1+\lambda}}&\mbox{if }0<\lambda\leq 2.
\end{array}\right.
\end{equation}}

\

{\bf Proof.} First note that if $x_1\to 0$ or $+\infty,$ then $x_2\to+\infty$
or $0$, and, in both cases, $\ds F_{\frac{1}{2}}(x_1,x_2)\to1.$ Now consider
the points $(x_1,x_2)\in H_\lambda$ for which
$\ds\frac{x_1}{(1+x_1)^{\frac{3}{2}}}=\frac{x_2}{(1+x_2)^{\frac{3}{2}}}.$
If $x_1=x_2,$ then $\ds F_{\frac{1}{2}}(x_1,x_2)=\frac{2}{\sqrt{1+\lambda}}.$
If $x_1\neq x_2,$ then $x_1+x_2=\lambda^2(\lambda^2-3)$ and since
$x_1x_2=\lambda^2$ we have $\lambda >2$ and
$\ds F_{\frac{1}{2}}(x_1,x_2)=\frac{\lambda}{\sqrt{\lambda^2-1}}.$ 
Hence (9) follows from the inequalities
$\ds\frac{\lambda}{\sqrt{\lambda^2-1}}>1$ and
$\ds\frac{\lambda}{\sqrt{\lambda^2-1}}\ge\frac{2}{\sqrt{1+\lambda}}$ for
$\lambda>2,$ and $\ds\frac{2}{\sqrt{1+\lambda}}>1$ for $\lambda\leq 2.$\qed

\

The dual analog of Proposition 2 is the following

\

{\bf Proposition 5.} {\it For any $\alpha\in (0,1)$ and any integer
$\ds n\ge\frac{1}{\alpha}+1$ we have
$$\sup_{H_\lambda}F_\alpha=\max(n-1,\frac{n}{(1+\lambda)^\alpha}).$$}

\

{\bf Proof.} Proceedings as in the proof of Proposition 2 it is enough to show
that $\ds(1+\frac{u}{\beta})^\beta<\frac{1-u^{\beta+1}}{1-u}$ for arbitrary
$u\in(0,1)$ and $\beta>1.$ This can be done in the same way as the proof of
the inequality (8).\hfill$\Box$

\

Finally, note that in the case $\alpha<0$ obviously the supremum of
$F_\alpha$ is equal to $+\infty.$

\

{\bf Remark 2.} As Example 2 suggests, if $\alpha\in(0,1)$ and
$\ds n<\frac{1}{\alpha}+1,$ then a result similar to Proposition 5 is not
true. The authors do not know the value of $\ds\sup_{H_\lambda}F_\alpha$
for such $\alpha$ and $n$.

\end{document}